\documentclass[12pt]{article}

\usepackage{amssymb, amsfonts, amsmath}
\usepackage{mathrsfs}
\usepackage{graphicx}
\usepackage[usenames,dvipsnames]{xcolor}
\usepackage{tikz}
\usepackage{float}
\usepackage{enumitem}
\usepackage{comment}
\usepackage{setspace}
\usepackage{caption}
\usepackage{subcaption}
\usepackage[hidelinks]{hyperref}
\usepackage{array}%
\usepackage{multicol}

\newcolumntype{C}[1]{>{\centering\let\newline\\\arraybackslash\hspace{0pt}}m{#1}}
\providecommand{\U}[1]{\protect\rule{.1in}{.1in}}

\newcommand{\Mod}[1]{\ (\text{mod}\ #1)}
\newcommand{\ig}[1]{\mathscr{I}(#1)}

\newcommand{\Dia}{\mathfrak{D}}

\newcommand{\wlat}{\mathfrak{L}}
\newcommand{\bcl}{\mathfrak{B}}

\newcommand{\cid}[1]{\left\langle #1 \right\rangle }
\newcommand{\bk}[1]{\left\lbrace  #1 \right\rbrace  }

\def\ieInner(#1,#2,#3,#4){#1 \overset{#2 #3}{\sim} #4}
\def\edge#1{\ieInner(#1)}
\def\ieInnerA(#1,#2,#3){\overset{#1 #2}{\sim} #3}

\def\sedge(#1,#2){#1\sim #2}
\def\ieInnerB(#1,#2,#3,#4,#5){#1 \overset{#2 #3}{\sim_{#5}} #4} 

\def\ieInnerC(#1,#2,#3){#1 {\sim_{#3}} #2 }

\def\ieInnerD(#1,#2,#3,#4){#1 \sim #2 \sim \dots \sim #3 \sim #4}

\def\ieInnerE(#1,#2){#1 \sim #2}

\def\ieInnerF(#1,#2){\overset{#1 #2}{\sim}}

\definecolor{ltsky}{RGB}{0,191,255}
\definecolor{ltteal}{RGB}{0, 128, 128 }
\definecolor{medOrch}{RGB}{122,55,139}
\definecolor{royalBlue}{RGB}{65,105,225}
\definecolor{forGreen}{RGB}{34,139,34}
\definecolor{dand}{RGB}{255,193,37}
\definecolor{lightBlue}{RGB}{176,226,255}
\definecolor{lgrey}{RGB}{209,209,209}
\definecolor{lgray}{gray}{0.95}
\definecolor{mgray}{gray}{0.40}
\definecolor{mmgray}{gray}{0.60}
\definecolor{zelim}{RGB}{7,163,82}
\definecolor{zelim2}{RGB}{5,114,57}

\usetikzlibrary{fit,positioning,calc, arrows, shapes}
\usetikzlibrary{decorations.pathreplacing}
\usetikzlibrary{backgrounds}
\tikzstyle{std}=[ circle, draw=black,fill=black,thick, inner sep=2pt, minimum size=2.5mm]
\tikzstyle{wstd}=[ circle, draw=black,fill=white,thick, inner sep=2pt, minimum size=2.5mm]
\tikzstyle{ir}=[ circle, draw=black,fill=green,thick,  inner sep=2pt, minimum size=2mm]
\tikzstyle{mp}=[circle, draw=black,fill=Dandelion,thick,  inner sep=2pt, minimum size=2mm]
\tikzstyle{bred}=[circle, draw=black,fill=red,thick,  inner sep=2pt, minimum size=2mm]
\tikzstyle{byellow}=[circle, draw=black,fill=yellow,very thick,  inner sep=2pt, minimum size=3mm]
\tikzstyle{bgray}=[circle, draw=black,fill=mmgray,thick,  inner sep=2pt, minimum size=2mm]
\tikzstyle{bzed}=[circle, draw=black,fill=zelim,thick,  inner sep=2pt, minimum size=2mm]
\tikzstyle{smred}=[ circle, draw=black,fill=red,thick,  inner sep=1pt, minimum size=1.5mm]
\tikzstyle{bblue}=[ circle, draw=black,fill=blue,thick,  inner sep=2pt, minimum size=2.5mm]
\tikzstyle{bgreen}=[ circle, draw=black,fill=green,thick,  inner sep=2pt, minimum size=2.5mm]
\tikzstyle{regRed}=[ circle, draw=black,fill=red,thick,  inner sep=2pt, minimum size=2.5mm]
\tikzstyle{sqRed}=[rectangle, draw=black,fill=red,thick,  inner sep=2pt, minimum size=2.5mm]
\tikzstyle{regG}=[ circle, draw=black,fill=green,thick,  inner sep=2pt, minimum size=2mm]
\tikzstyle{dred}=[diamond, draw=black,fill=red,thick,  inner sep=2pt, minimum size=2mm]
\tikzstyle{tr}=[color=black, style=dotted]
\tikzstyle{sp}=[color=ProcessBlue, line width = 2pt]

\tikzstyle{bline}=[color=blue, line width = 1pt]
\tikzstyle{rline}=[color=red, line width = 1pt]
\tikzstyle{rlinelt}=[color=red, line width = 1pt,opacity=0.2]
\tikzstyle{rlinemb}=[color=red, line width = 1pt, densely dashdotted]
\tikzstyle{gline}=[color=mmgray, line width = 1pt]
\tikzstyle{zline}=[color=zelim, line width = 1pt]
\tikzstyle{ll}=[color=gray]

\tikzstyle{vertex}=[circle, fill=black, inner sep= 0, minimum size = 4]

\tikzset{broadArrow/.style={single arrow, fill=red!50, anchor=base, align=center,text width=2.8cm}}

\pgfdeclarelayer{blob}    
\pgfdeclarelayer{bedge}    
\pgfsetlayers{bedge,blob,main}

\usepackage{authblk}

\newtheorem{theorem}{Theorem}[section]

\newtheorem{lemma}[theorem]{Lemma}

\newtheorem{obs}[theorem]{Observation}

\newtheorem{problem}{Problem}

\newenvironment{proof}[1][Proof]{\noindent\textbf{#1.} }{\ \hfill \rule{0.5em}{0.5em}}

\begin{document}

\title{\textbf{The $i$-Graphs of Paths and Cycles}}


\author{R. C. Brewster\thanks{Funded by Discovery Grants from the Natural Sciences and
Engineering Research Council of Canada, RGPIN-2014-04760, RGPIN-03930-2020.}}
\affil{Department of Mathematics and Statistics\\Thompson Rivers University\\805 TRU Way\\Kamloops, B.C.\\ \textsc{Canada} V2C 0C8}
\author{C. M. Mynhardt$^{*}$}
\author{L. E. Teshima}
\affil{Department of Mathematics and Statistics\\University of Victoria\\PO BOX 1700 STN CSC\\Victoria, B.C.\\\textsc{Canada} V8W 2Y2 \authorcr {\small rbrewster@tru.ca, kieka@uvic.ca, lteshima@uvic.ca}}

\maketitle

\begin{abstract}

The independent domination number $i(G)$ of a graph $G$ is the minimum cardinality of a maximal independent set 
of $G$, also 
called an $i(G)$-set.  The $i$-graph of $G$, denoted $\ig{G}$, is the graph whose vertices correspond to the 
$i(G)$-sets, and where two $i(G)$-sets are adjacent if and only if they differ by two adjacent vertices. Although
not all graphs are $i$-graph realizable, that is, given a target graph $H$, there does not necessarily exist a source graph $G$ such that $H \cong \ig{G}$, all graphs have $i$-graphs.  We determine the $i$-graphs of paths and cycles and, in the case of cycles, discuss the Hamiltonicity of these $i$-graphs.
\end{abstract}

\noindent\textbf{Keywords:\hspace{0.1in}}independent domination number, graph reconfiguration, $i$-graph, $i$-graph of a path, $i$-graph  of a cycle

\noindent\textbf{AMS Subject Classification Number 2020:\hspace{0.1in}}05C69	

\section{Introduction}
\label{sec:intro}
Paths and cycles are among the simplest graph classes and one can easily determine their various domination-type numbers. We consider their minimum independent dominating sets to be the vertices of a new graph, called an ``$i$-graph'', in which two vertices are adjacent whenever the symmetric difference of their corresponding sets consists of two vertices that are adjacent in the original graph. We also precisely determine those cycles whose $i$-graphs are Hamiltonian or traceable, i.e., have Hamiltonian paths. 

The $i$-graph of a graph is an example of a ``reconfiguration graph''. In graph theory, reconfiguration problems are often concerned with solutions to a specific problem that are
vertex/edge subsets of a graph. When this is the case, the reconfiguration problem can be viewed as 
a token manipulation problem, where
a solution subset is represented by placing a token at each vertex or edge of the subset.
The reconfiguration step for vertex subsets can be of one of three variants (edge
subsets are handled analogously):

\begin{itemize}
\item[$\vartriangleright $] \textbf{Token Slide (TS) Model}: A single token
is slid along an edge between adjacent vertices.

\item[$\vartriangleright $] \textbf{Token Jump (TJ) Model}: A single token
jumps from one vertex to another (without the vertices necessarily being
adjacent).

\item[$\vartriangleright $] \textbf{Token Addition/Removal (TAR) Model:} A
single token can either be added to a vertex or be removed from a vertex.
\end{itemize}

To represent the many possible solutions in a reconfiguration problem, each solution can
be represented as a vertex of a new graph, referred to as a \textit{reconfiguration graph}, where adjacency between vertices follows one of the three token adjacency models, producing the
\textit{slide graph}, the \textit{jump graph}, or the \textit{TAR graph}, respectively. Here we consider the token slide model reconfiguration step applied to the minimum independent dominating sets of paths and cycles.
We refer the reader to \cite{MN20} for a survey on reconfiguration of colourings and dominating sets in graphs. 
See also the survey~\cite{Nishimura18} for results on independent set reconfiguration.

We use the standard notation of $\alpha(G)$ for the independence number and $\gamma (G)$ for the domination number 
(cardinality of a minimum dominating set) of a graph $G$. 
The \emph{independent domination number} $i(G)$ of $G$ is the
minimum cardinality of a maximal independent set of $G$, or, equivalently,
the minimum cardinality of an independent domination set of $G$. An independent dominating set of $G$ of cardinality $i(G)$ is also called an $i$-\emph{set} of $G$, or an $i(G)$-\emph{set}. In general, we follow the notation of \cite{CLZ}. 

\subsection*{$i$-Graphs}
The \emph{$i$-graph} of a graph $G$, denoted $\ig{G}=(V(\ig{G}),E(\ig{G}))$, is the graph with vertices 
representing the minimum independent dominating sets of $G$ (that is, the {$i$-sets} of $G$), and where $u,v\in V(\ig{G})$, corresponding to the $i(G)$-sets  $S_u$ and $S_v$, respectively, are 
adjacent in $\ig{G}$ if and only if there exists  $xy \in E(G)$ such that $S_u = (S_v-x)\cup\{y\}$.  That is, adjacency in $\ig{G}$ follows a token slide model. 

We say $H$ \emph{is an $i$-graph}, or is $i$-\emph{graph realizable}, if there exists some graph $G$ such that 
$\ig{G}\cong H$. Moreover, we refer to $G$ as the \emph{seed graph} of the $i$-graph $H$.  
Going forward, we mildly abuse notation to denote both the $i$-set $X$ of $G$ and its corresponding 
vertex in $H$ as $X$, so that $X \subseteq V(G)$ and $X \in V(H)$. 

We observe that only the token slide model is relevant for $i$-set reconfiguration. Imagine that there is a token on each vertex of an $i$-set $S$ of $G$.  Then $S$ is adjacent, in $\ig{G}$, to an $i(G)$-set $S'$ if and only if a single token can be slid along an edge of $G$ to transform $S$ into $S'$.  
A token is said to be \emph{frozen} (in any reconfiguration model) if there are no available vertices to which it can slide.

In acknowledgment of the slide-action in $i$-graphs, given $i$-sets  $X = \{x_1,x_2,\dots,x_k\}$ and 
$Y=\{y_1,x_2,\dots x_k\}$ of $G$ with $x_1y_1 \in E(G)$, we denote the adjacency of $X$ and $Y$ in $\ig{G}$ 
as $\edge{X,x_1,y_1,Y}$, where we imagine transforming the $i$-set $X$ into $Y$ by sliding the token at $x_1$ 
along an edge to $y_1$.  
More generally, we use $x \sim y$ to denote the adjacency of vertices $x$ and $y$ (and $x\not \sim y$ to 
denote non-adjacency); this is used in the context of both the seed graph and the target graph.

The study of $i$-graphs was initiated by L. E. Teshima in \cite{LauraD}. In the paper \cite{BMT1} based on this
work, the authors investigated $i$-graph realizability and proved a number of results concerning the adjacency of
vertices in an $i$-graph and the structure of their associated $i$-sets in the seed graph. They presented the 
three smallest graphs that are not $i$-graphs, namely the diamond graph $\Dia = K_4-e$, $K_{2,3}$ and the graph
$\kappa$, which is $K_{2,3}$ with an edge subdivided. They showed that several common graph classes, like trees
and cycles, are $i$-graphs. They demonstrated that known $i$-graphs can be used to construct new $i$-graphs and 
applied these results to build other classes of $i$-graphs, such as block graphs, hypercubes, forests, 
cacti, and unicyclic graphs.

The diamond $\Dia$, $K_{2,3}$, and $\kappa$, mentioned above, are examples of \emph{theta graphs}: graphs that 
are the union of three internally disjoint nontrivial paths with the same two distinct end vertices.  
The problem of characterizing theta graphs that are $i$-graph realizable was fully resolved in \cite{LauraD} 
and also reported in  \cite{BMT3}.

Here we consider the opposing question: given a graph $G$, what is the structure of $\ig{G}$? 
The exact structure of the resulting $i$-graph can vary among families of graphs from the simplest 
isolated vertex to surprisingly complex structures.  
We illustrate this statement by examining the $i$-graphs of two of the most famous classes of graphs: 
paths and cycles. We count the number of distinct $i$-sets of the path $P_n$ and the cycle $C_n$ in Sections \ref{sec:i:numbPaths} 
and \ref{sec:i:numCycles}, respectively.  
That is, we determine $ \left|  V \left(\ig{P_n} \right)  \right|$ for $n\geq 1$  and $\left|  V \left(\ig{C_n} \right)  \right|$ for $n \geq 3$.  In Sections \ref{sec:i:Pn} and \ref{sec:i:cycles} we determine the $i$-graphs of paths and cycles, respectively. Then, in Section \ref{sec:i:hamCycle}, we resolve the question of which cycles have Hamiltonian $i$-graphs, and for those that do not, we determine, in Section \ref{subsec:i:C3k1Trac}, which have \emph{traceable} $i$-graphs, that is, have $i$-graphs that admit Hamiltonian paths.  

We note that a description of the $i$-graphs for $P_n$ and $C_n$ appears in \cite{STA19}.  Our development uses a different description of the $i$-sets which gives shorter proofs and allows use to establish results on the Hamiltonicity of the $i$-graphs of $C_n$.  Results on the Hamiltonicity of other domination reconfiguration problems
appears in~\cite{Adaricheva2021}.

\section{The $i$-Graphs of Paths} 
\label{sec:paths}

We assume that the vertices of the path $P_{n}$  are labelled as 
$P_{n} = (v_1,v_2,\dots, v_{n})$ in the natural order. 
Given that we are discussing $i$-sets, which are both independent and dominating, if $X$ is an $i$-set of $P_n$,
then two consecutive vertices of $X$ are separated by one or two vertices of $P_n-X$; the different interval 
lengths between these consecutive vertices of $X$ therefore correspond to the different $i$-sets of $P_n$.  
This provides our method for counting the distinct $i$-sets of $P_n$. 

\subsection{The Number of $i$-Sets of Paths} \label{sec:i:numbPaths}

To begin, recall the following well-known result regarding the  independent domination number for both paths and cycles.

\begin{lemma} \label{lem:i:iPC} \emph{\cite{GH13}}
	For the path and cycle, $i(P_n) = i(C_n) = \left \lceil n/3 \right \rceil $.
\end{lemma}

\noindent 
For the path $P_n$, let $t = i(P_n)+1$.  Given an arbitrary $i$-set $X$ of $P_n$, $V(P_n)-X$ is partitioned into $t$ intervals $X_1, X_2, \dots, X_t$. By our above remarks, for $2 \leq i \leq t-1$, $x_i := |X_i|$ is $1$ or $2$.  We call $X_i$ \emph{small} if $x_i = 1$ and \emph{large} if $x_i = 2$.  Similarly $X_1$ and $X_t$ are small and large when they have size $0$ and $1$ respectively.
An example for $P_{10}$ with sets $X_1, X_2, \dots, X_5$ is given below in Figure~\ref{fig:i:GF}.  In particular, notice that $X_1 = \varnothing$. Both $X_1$ and $X_2$ are small while the other intervals are large.

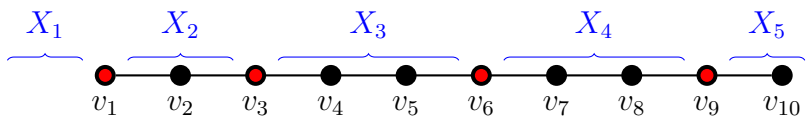
\begin{figure}[H] \centering	
	\begin{tikzpicture}			 		
		
		\tikzset{
			position label/.style={
				above = 1pt,
				text height = 1ex,
				text depth = 1ex
			},
			brace/.style={
				decoration={brace, mirror},
				decorate
			}
		}

		\coordinate (cent) at (0,0) {};
		\coordinate (v0) at (0,0) {};
		
		\foreach \i / \j in {1,2,3,4,5,6,7,8,9,10}
		{
			\path(cent) ++(0: \i*10 mm) node[std,label={ 270:$v_{\i}$}]  (v\i) {};
		}
		\draw[thick] (v1)--(v10);

		\foreach \i / \j in {1,3,6,9}
		\node[bred] (w\i) at (v\i) {};

		\foreach \i / \j in {0,1,2,3,4,5,6,7,8,9,10}
		{
			\path(v\i) ++(-3mm,2mm) coordinate (l\i) {};
			\path(v\i) ++(3mm,2mm) coordinate (r\i) {};
		}
		
		\draw [decoration={brace}, decorate, color=blue] (l0.north) -- node [position label, pos=0.5] {$X_1$} (l1.north);
		\draw [decoration={brace}, decorate, color=blue] (r1.north) -- node [position label, pos=0.5] {$X_2$} (l3.north);
		\draw [decoration={brace}, decorate, color=blue] (r3.north) -- node [position label, pos=0.5] {$X_3$} (l6.north);
		\draw [decoration={brace}, decorate, color=blue] (r6.north) -- node [position label, pos=0.5] {$X_4$} (l9.north);
		\draw [decoration={brace}, decorate, color=blue] (r9.north) -- node [position label, pos=0.5] {$X_5$} (r10.north);		
	\end{tikzpicture}   
	
	\caption{The sets $X_1,X_2,\dots,X_5$ of $P_{10}$.}
	\label{fig:i:GF}	
\end{figure}	

As each vertex of $P_n$ belongs to $X$ or to some $X_i$, and the corresponding $x_i$ is bounded above by $1$ ($i=1$ or $t$) or by $2$ ($2 \leq i \leq t-2$), 
we have the following inequality.
$$
n=|V(P_n)| = |X| + \sum_{i=1}^t x_i \leq (t-1) + 1 + 2(t-2) + 1 = 3t-3.
$$
Since $\lceil n/3 \rceil = t-1$, it is straightforward to see that
$$
n = \left\{ \begin{array}{ll}
3t-3 & \mbox{ if } n = 3k \\
3t-5 & \mbox{ if } n = 3k+1 \\
3t-4 & \mbox{ if } n = 3k+2.
\end{array} \right.
.
$$
From this we can conclude that for $n=3k$, each $X_i$ is large and $P_n$ has a unique $i$-set.  When $n=3k+1$, there are exactly two sets $X_i$ and $X_j$ that are small and $P_n$ has exactly ${t \choose 2} = {{k+2} \choose 2}$ $i$-sets.  Finally when $n=3k+2$, there is one small $X_i$ and $P_n$ has exactly $t=k+2$ $i$-sets. In summary,

\begin{lemma} \label{lem:i:pathSize}
	For $n \geq 1$, the order of $\ig{P_n}$ is 	
	\begin{align*}
		\left|  V \left(\ig{P_n} \right)  \right|  = 
		\begin{cases}
			1 & \text{ if } \ n = 3k\\ 
			\binom{k+2}{2} & \text{ if }  \ n = 3k+1 \\
			k+2 & \text{ if }  \ n=3k+2.
		\end{cases}
	\end{align*}
\end{lemma}

\subsection{The $i$-Graph of $P_n$} \label{sec:i:Pn}

To see the structure of $\ig{P_n}$, notice that the $i$-set tokens on $P_n$ that are free to slide are very limited. For example, in Figure \ref{fig:i:P10}	below, we have two different $i$-sets on $P_{10}$.  In the first case, the token at $v_6$ is frozen as $X_3$ and $X_4$ are both large.  Moving the token would create an interval exceeding the large size and thus leave an undominated vertex.  The token at $v_3$ can move to the right as $X_2$ is small and $X_3$ is large.  After such a move, $X_2$ is large and $X_3$ is small.

\begin{figure}[H] \centering	
	\begin{tikzpicture}			
		\coordinate (cent) at (0,0) {};
		
		\foreach \i / \j in {1,2,3,4,5,6,7,8,9,10}
		{
			\path(cent) ++(0: \i*8 mm) node[std,,label={ 270:$v_{\i}$}]  (v\i) {};
		}
		\draw[thick] (v1)--(v10);

		\foreach \i / \j in {1,3,6,9}
		\node[bred] (w\i) at (v\i) {};
		
		\draw[ltteal,thick, ->]  (v3) to [out=45,in=135]  (v4);

	\end{tikzpicture} 
	
	\vspace{1cm}
	
	\begin{tikzpicture}			
		\coordinate (cent) at (0,0) {};
		
		\foreach \i / \j in {1,2,3,4,5,6,7,8,9,10}
		{
			\path(cent) ++(0: \i*8 mm) node[std,,label={ 270:$v_{\i}$}]  (v\i) {};
		}
		\draw[thick] (v1)--(v10);

		\foreach \i / \j in {2,4,7,9}
		\node[bred] (w\i) at (v\i) {};
		
		\draw[ltteal,thick, <-]  (v1) to [out=45,in=135]  (v2);
		\draw[ltteal,thick, ->]  (v4) to [out=45,in=135]  (v5);
		\draw[ltteal,thick, <-]  (v6) to [out=45,in=135]  (v7);
		\draw[ltteal,thick, ->]  (v9) to [out=45,in=135]  (v10);
		
	\end{tikzpicture} 
	
	\caption{Two $i$-sets of $P_{10}$. In the first, the $i$-set $(1,2)$ becomes $(1,3)$ after the token slide.  In the second, the allowed token slides show the $i$-set $(2,4)$ is adjacent to $(1,4), (3,4), (2,3),$ and $(2,5)$ in $\ig{P_{10}}$.}
	\label{fig:i:P10}
\end{figure}
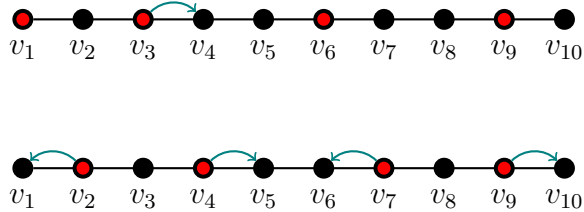	

Our key observation is that a token between $X_i$ and $X_{i+1}$ can slide to the right if and only if $X_i$ is small and $X_{i+1}$ is large, after which $X_i$ is large and $X_{i+1}$ is small, and vice versa for sliding to the left.

When $n=3k+2$, each $i$-set has a single small interval which uniquely identifies the $i$-set.  From our observation we conclude two $i$-sets are adjacent in $\ig{P_n}$ if the indices of their small intervals differ by one.  Thus, $\ig{P_n} = P_{k+2}$.
When $n=3k+1$, each $i$-set has two small intervals, say $X_i$ and $X_j$, and thus the pair $(i,j)$ uniquely identifies the $i$-set, for $1 \leq i < j \leq k+2$.  Let $\wlat_k$ denote the subgraph of the integer lattice induced by $\{ (i, j):1 \leq i < j \leq k+2 \}$.  The pair $(i, j)$ is adjacent to $(i',j')$ when $1 \leq i' < j' \leq k+2$ and exactly one of $|i-i'|=1$ or $|j-j'|=1$ holds. From our observation, we conclude that $\ig{P_n} = \wlat_k$.  
These results are summarized below.

\begin{theorem} \label{thm:i:Pn}
  The $i$-graph for the path $P_n$ is given by
	\begin{align*}
		\ig{P_n} = 	
		\begin{cases}
			K_1 & \text{ if } n=3k, \\
			\wlat_k & \text{ if } n=3k+1, \\
			P_{k+2} & \text{ if } n=3k+2.
		\end{cases}
	\end{align*}
\end{theorem}

An example of the $i$-graph $\wlat_3$ is shown in Figure~\ref{fig:wlat}.
\begin{figure}[H] \centering	
	\begin{tikzpicture}[scale=1.6]	
      \foreach \i in {2,3,4,5}
        \node[std,label={180:$(1,\i)$}] (a\i) at (1,\i) {};
      \foreach \i in {3,4,5}
        \node[std,label={315:$(2,\i)$}] (b\i) at (2,\i) {};  
      \foreach \i in {4,5}
        \node[std,label={315:$(3,\i)$}] (c\i) at (3,\i) {};        
      \node[std,label={315:$(4,5)$}] (d5) at (4,5) {};        
    \draw[thick] (a2)--(a5) (b3)--(b5) (c4)--(c5) (a3)--(b3)
        (a4)--(c4) (a5)--(d5);
    \end{tikzpicture}
\caption{The graph $\ig{P_{10}} = \wlat_3$.}\label{fig:wlat}  
\end{figure}    

\section{The $i$-Graphs of Cycles}

Contrary to our conventions in the previous sections on paths, we assume that all cycles have labelled vertex set $V(C_n) = (v_0,v_1,\dots,v_{n-1})$ (according to some orientation of the cycle).

\subsection{The Number of $i$-Sets of Cycles}	\label{sec:i:numCycles}

Our method for counting the $i$-sets of cycles is similar to our approach to paths.  An $i$-set $X$ of $C_n$ again partitions $V(C_n)-X$ into sets $X_1, X_2, \dots, X_t$. Letting $x_j = |X_j|$, we have $1 \leq x_j \leq 2$ with exactly $0, 1,$ or $2$ of the $x_j = 1$ when $n = 3k, 3k+2,$ or $3k+1$ respectively.  The rotational symmetry of $C_n$ however complicates the counting. For example, when $n=3k$ each $x_j = 2$ uniquely determining the gaps between vertices of $X$, but there are 3 $i$-sets of $C_n$.  




Removing the edge $e=v_0v_{n-1}$ of $C_n$ gives the path $P_n = C_n - e = (v_0, \dots, v_{n-1})$.  Any $i$-set $X$ of $C_n$ corresponds to either an $i$-set of $P_n$ or to a ``near $i$-set of $P_n$'', in which 

\begin{enumerate}[label=(\roman*)]
	\item $v_2$ is the first vertex and $v_{n-1}$ is the last vertex of $X$ on $P_n$, or
	\item $v_0$ is the first vertex and $v_{n-3}$ is the last vertex of $X$ on $P_n$.
\end{enumerate}

These three varieties of $i$-sets are pairwise disjoint; therefore, by counting the number of $i$-sets and near $i$-sets of $P_n$ we obtain the number of $i$-sets of $C_n$.  Recall $t = i(P_n)+1$.  Define $r = \sum_{j=1}^t x_j = n-i(P_n)$.
Thus, the number of $i$-sets of $C_n$ equals the sum of the number of integer solutions to the following equations:
\begin{align} \label{eq:i:cycleEQ}
	x_1 + \dots + x_t = r, & \text{ where } x_1=0 \text{ and } 1 \leq x_j \leq 2 \text{ for } j \in \{ 2,\dots,t \},\\
	x_1 + \dots + x_t = r, & \text{ where } x_t=0 \text{ and } 1 \leq x_j \leq 2 \text{ for } j \in \{1,\dots,t-1\},\\
	x_1 + \dots + x_t = r, & \text{ where } x_1=x_t=1 \text{ and } 1 \leq x_j \leq 2 \text{ for } j \in \{2,\dots,t-1\}.
\end{align}

The corresponding generating function is 
\begin{align} \label{eq:i:cycleGF}
	(x + x^2)^{t-1} + (x+x^2)^{t-1} + x^2(x+x^2)^{t-2} = 2x^{t-1}(1+x)^{t-1} + x^t(1+x)^{t-2}.
\end{align}

\begin{enumerate}[label=$\bullet$]
	\item If $n=3k$, then $t=k+1$ and $r=2k$. Thus, we require the coefficient of $x^{2k}$ in $2x^{k}(1+x)^{k} + x^{k+1}(1+x)^{k-1}$, which is 3.
 
    \item If $n=3k+1$, then $t = i(P_n) + 1 = k+2$ and $r=2k$.  Thus, we require the coefficient of $x^{2k}$ in (\ref{eq:i:cycleGF}),  that is, in $2x^{k+1}(1+x)^{k+1} + x^{k+2}(1+x)^k$, which is $2 \binom{k+1}{k-1} + \binom{k}{k-2} = k(3k+1)/2$.
	
	\item If $n = 3k+2$, then $t=k+2$ and $r=2k+1$.  Thus, we require the coefficient of $x^{2k+1}$ in $2x^{k+1}(1+x)^{k+1}+x^{k+2}(1+x)^k$, which is $2\binom{k+1}{k} + \binom{k}{k-1} = 3k+2 = n$.
	
	
\end{enumerate}

We summarize these results in the lemma below.

\begin{lemma} \label{lem:i:cycleSize}
	For $n \geq 3$, the order of $\ig{C_n}$ is 
	
	\begin{align*}
		\left|  V \left(\ig{C_n} \right)  \right|  = 
		\begin{cases}
			3 & \text{ if } \ n = 3k \\ 
			k(3k+1)/2 & \text{ if }  \ n = 3k+1 \\
			n & \text{ if }  \ n=3k+2.
		\end{cases}
	\end{align*}
\end{lemma}

\subsection{The $i$-Graphs of Cycles}	\label{sec:i:cycles}


Immediately, Lemma \ref{lem:i:cycleSize} shows that some of the $i$-graphs for $C_n$ are fairly straight-forward.  For $C_{3}$, this is a complete graph, and hence $\ig{C_3} = \ig{K_3} = C_3$.  When $k \geq 2$, $C_{3k}$ has three distinct $i$-sets. In each case, each $i$-set vertex has two non-$i$-set vertices between it and the next $i$-set vertex. Thus each $i$-set vertex is frozen and $\ig{C_n}$ consists of three singletons.   

For $C_{n}$ with $n \equiv 2 \Mod{3}$, say $n=3k+2$, each $i$-set of $C_n$ contains exactly one pair of vertices $v_{j-1}$ and $v_{j+1}$ that are separated by exactly one vertex, $v_j$, not in the $i$-set (the common neighbour of these two vertices), while all other pairs  of consecutive $i$-set vertices are separated by exactly two vertices not in the $i$-set.  Hence, the $i$-set has exactly two vertices, namely $v_{j-1}$ and $v_{j+1}$, that are not frozen, and each of them can slide in only one direction.  The vertex $v_{j-1}$ can move to $v_{j-2}$, and $v_{j+1}$ can move to $v_{j+2}$.  In the former case,  the vertex that has two neighbours  in the $i$-set is now $v_{j-3}$, while in the latter, it is $v_{j+3}$.  As a result, $\ig{C_n}$ is 2-regular and since 3 is coprime to $n$, the sequence $v_j$, $v_{j+3}$, $v_{j+6}$, $\dots$, will visit each vertex in $C_n$  after $n$-slides.  Thus, each $i$-set is generated. Hence, $\ig{C_n}$ is 2-regular and connected; we conclude it is the cycle $C_n$.

We provide examples of $C_5$ and $C_8$ in Figures \ref{fig:i:C5} and \ref{fig:i:C8} below.

\begin{figure}[H] \centering
	\begin{subfigure}{.4\textwidth}
		\begin{tikzpicture}	[scale=.8]		
			\coordinate (cent) at (0,0) {};

			\foreach \i  in {0,1,2,3,4} 
			{
				\node[std,label={90-\i*72: ${v_{\i}}$}] (v\i) at (90-\i*72:2cm) {}; 
			}

			\draw[thick] (v0)--(v1)--(v2)--(v3)--(v4)--(v0);
			
			\foreach \i  in {1,3} 
			{
				\node[regRed] (r\i) at (v\i) {};
			}
		\end{tikzpicture} 
		\caption{$C_5$ with the $i$-set $\{v_1,v_3\}$ in red.}
		\label{subfig:i:C5}
	\end{subfigure}
	\hspace{10mm}
	\begin{subfigure}{.4\textwidth}
		\begin{tikzpicture}	[scale=.8]		
			\coordinate (cent) at (0,0) {};

			\foreach \i / \x / \y  in {0/0/2,1/0/3,2/1/3,3/1/4,4/2/4} 
			{
				\node[std,label={90-\i*72: $\{v_{\x},v_{\y}\}$}] (v\i) at (90-\i*72:2cm) {}; 
			}	
			
			\draw[thick] (v0)--(v1)--(v2)--(v3)--(v4)--(v0);
			
			\foreach \i  in {2} 
			{
				\node[regRed] (r\i) at (v\i) {};
			}
		\end{tikzpicture} 
		\caption{$\ig{C_5}$ with the vertex corresponding to the $i$-set $\{v_1,v_3\}$ in red. }
		\label{subfig:i:C5i}
	\end{subfigure}
	\caption{$C_5$ and its $i$-graph.}
	\label{fig:i:C5}
\end{figure}
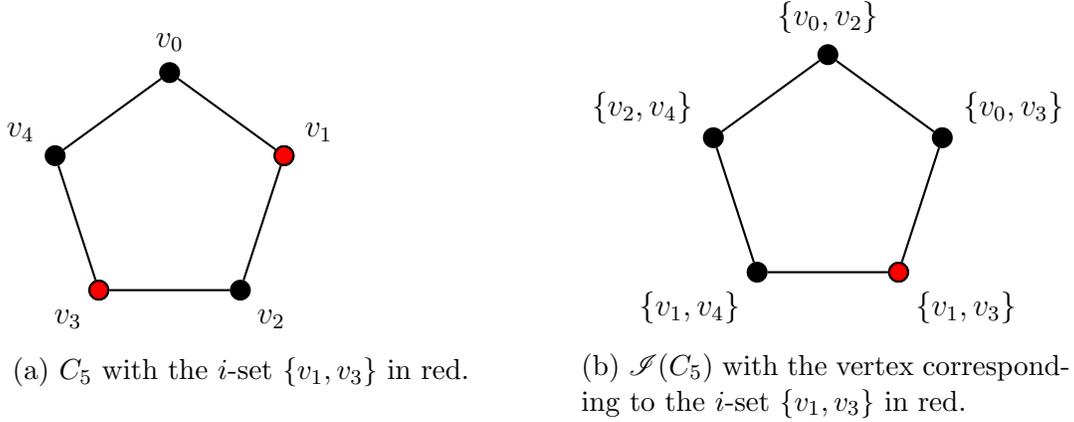

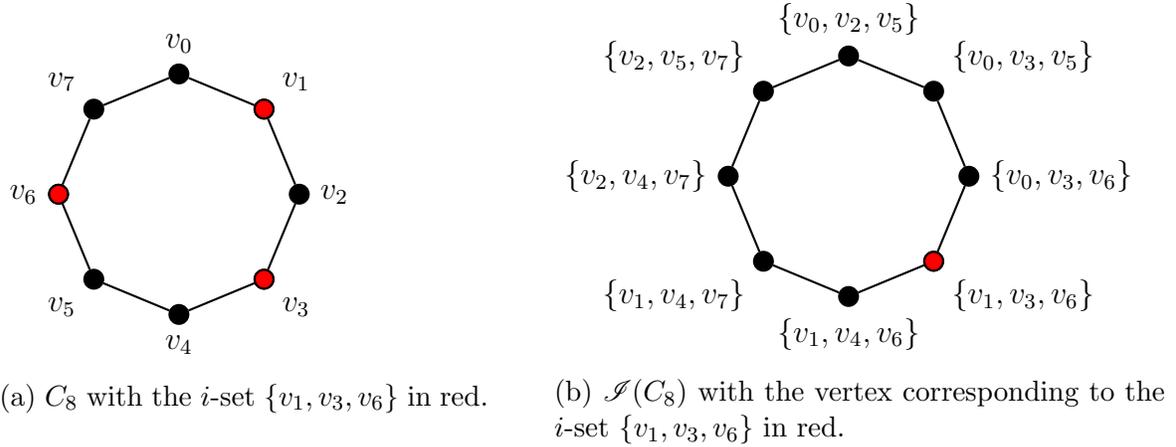
\begin{figure}[H] \centering
	\begin{subfigure}{.4\textwidth}
		\begin{tikzpicture}	[scale=.8]		
			\coordinate (cent) at (0,0) {};

			\foreach \i  in {0,1,2,3,4,5,6,7} 
			{
				\node[std,label={90-\i*45: ${v_{\i}}$}] (v\i) at (90-\i*45:2cm) {}; 
			}

			\draw[thick] (v0)--(v1)--(v2)--(v3)--(v4)--(v5)--(v6)--(v7)--(v0);
			
			\foreach \i  in {1,3,6} 
			{
				\node[regRed] (r\i) at (v\i) {};
			}
		\end{tikzpicture} 
		\caption{$C_8$ with the $i$-set $\{v_1,v_3,v_6\}$ in red.}
		\label{subfig:i:C8}
	\end{subfigure}
	\hspace{6mm}
	\begin{subfigure}{.5\textwidth}
		\begin{tikzpicture}	[scale=.8]		
			\coordinate (cent) at (0,0) {};

			\foreach \i / \x / \y / \z in {0/0/2/5,1/0/3/5,2/0/3/6,3/1/3/6,4/1/4/6,5/1/4/7,6/2/4/7,7/2/5/7} 
			{
				\node[std,label={90-\i*45: $\{v_{\x},v_{\y},v_{\z}\}$}] (v\i) at (90-\i*45:2cm) {}; 
			}	
			
			\draw[thick] (v0)--(v1)--(v2)--(v3)--(v4)--(v5)--(v6)--(v7)--(v0);
			
			\foreach \i  in {3} 
			{
				\node[regRed] (r\i) at (v\i) {};
			}
		\end{tikzpicture} 
		\caption{$\ig{C_8}$ with the vertex corresponding to the $i$-set $\{v_1,v_3,v_6\}$ in red. }
		\label{subfig:i:C8i}
	\end{subfigure}
	\caption{$C_8$ and its $i$-graph.}
	\label{fig:i:C8}
\end{figure}

We summarize these results in the lemma below.  
\begin{lemma} \label{lem:i:2Cn}
	For $n \geq 3$, $k \geq 0$,	
	\begin{align*}
		\ig{C_n} \cong 	
		\begin{cases}
			K_3 & \text{ if } n=3 \\
			3K_1 & \text{ if } n=3k \geq 6 \\
			C_n & \text{ if } n=3k+2.
		\end{cases}
	\end{align*}
\end{lemma}


Once again, the case that requires deeper analysis is $n = 3k+1$ for $k \geq 1$.  To help us tackle this final case, we first introduce a new notation for referencing the $i$-sets of cycles.  

For cycles $C_{3k+1}$, given any $i$-set $X$, there are exactly two vertices in $C_{3k+1} - X$ that are doubly dominated (that is, are adjacent to two different vertices of $X$).  Rather than referring to the $i$-set by its elements, which given a large $k$, could be numerous, we instead refer to the $i$-sets by these two unique vertices.  For clarity, we use a wide-angled bracketed notation when using this convention.  Figure \ref{fig:i:C13}	below illustrates this system:  for (a), rather than calling the $i$-set $X=\{v_1,v_3,v_6,v_9,v_{12}\}$, we refer to it as $\left\langle 0,2\right\rangle = \left\langle 2,0\right\rangle $.  Similarly in (b), instead of $Y=\{v_1,v_4,v_6,v_9,v_{12}\}$, we denote the $i$-set as $\left\langle 0,5\right\rangle $.

\begin{figure}[H] \centering	
	
	\begin{subfigure}{.35\textwidth}
		\begin{tikzpicture}	[scale=.8]		
			\coordinate (cent) at (0,0) {};

			\foreach \i  in {1,3,4,5,6,7,8,9,10,11,12} 
			{
				\node[std,label={90-\i*27.6923: ${\i}$}] (v\i) at (90-\i*27.6923:2cm) {}; 
			}

			\foreach \i  in {0,2} 
			{
				
				\node[std,label={90-\i*27.6923: $\color{blue} \fbox{{\i}}$}] (v\i) at (90-\i*27.6923:2cm) {}; 
				
			}	
			
			\draw[thick] (v0)--(v1)--(v2)--(v3)--(v4)--(v5)--(v6)--(v7)--(v8)--(v9)--(v10)--(v11)--(v12)--(v0);
			
			\foreach \i  in {1,3,6,9,12} 
			{
				\node[regRed] (r\i) at (v\i) {};
			}
		\end{tikzpicture} 
		\caption{$X=\{v_1,v_3,v_6,v_9,v_{12}\}$, \\ denoted as $\left\langle 0,2\right\rangle$.}
		\label{subfig:i:C13b}
	\end{subfigure}
	\hspace{3cm}
	\begin{subfigure}{.35\textwidth}
		\begin{tikzpicture}	[scale=.8]		
			\coordinate (cent) at (0,0) {};

			\foreach \i  in {1,2,3,4,6,7,8,9,10,11,12} 
			{
				\node[std,label={90-\i*27.6923: ${\i}$}] (v\i) at (90-\i*27.6923:2cm) {}; 
			}

			\foreach \i  in {0,5} 
			{
				
				\node[std,label={90-\i*27.6923: $\color{blue} \fbox{{\i}}$}] (v\i) at (90-\i*27.6923:2cm) {}; 
				
			}	
			
			\draw[thick] (v0)--(v1)--(v2)--(v3)--(v4)--(v5)--(v6)--(v7)--(v8)--(v9)--(v10)--(v11)--(v12)--(v0);
			
			\foreach \i  in {1,4,6,9,12} 
			{
				\node[regRed] (r\i) at (v\i) {};
			}
		\end{tikzpicture} 
		\caption{$Y=\{v_1,v_4,v_6,v_9,v_{12}\}$, \\ denoted as $\left\langle 0,5\right\rangle$.}
		\label{subfig:i:C13a}
	\end{subfigure}
	\caption{Two $i$-sets of $C_{13}$.}
	\label{fig:i:C13}
\end{figure}
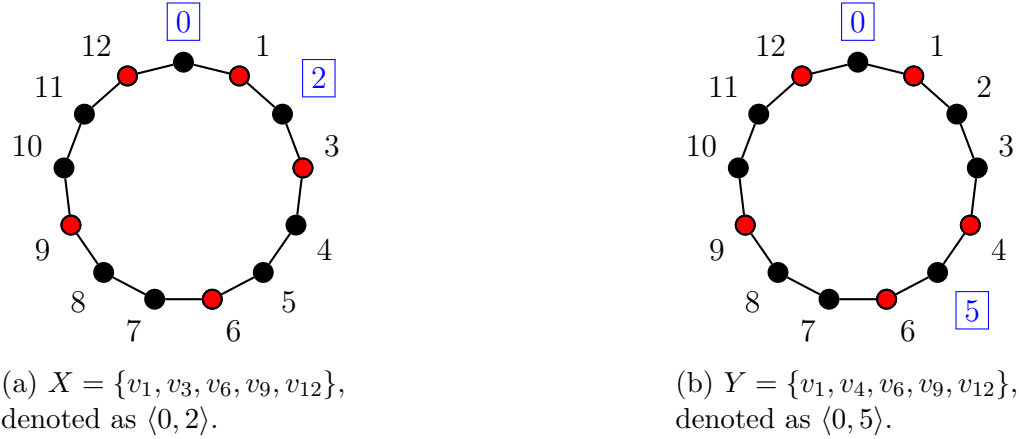	

With this new labelling in place, we move now to a proposed family of graphs that we have dubbed \emph{bracelet graphs}, $\bcl_{k}$.
The vertex set of $\bcl_{k}$ consists of all distinct $2$-subsets
$\{j,\ell\}$ of $\{0,1,\dots,3k\}$ such that $0\leq j\leq3k$ and $\ell\equiv
j+3s+2(\operatorname{mod}\ 3k+1)$, $s\in\{0,1,\dots,k-1\}$. For example, the
subsets containing $0$ are $\{0,2\},\{0,5\},\dots,\{0,3k-4\},\{0,3k-1\}$.

To simplify notation, assume the vertices of $C_{3k+1}$ are labelled $0,1,\dots,3k$ in clockwise
order, as illustrated in Figure \ref{fig:i:C13}. In $\mathfrak{B}_{k}$, the neighbours of
the vertex $\{j,\ell\}$ are described below.
\begin{enumerate}
	\item Suppose $j-\ell\equiv2$ or $-2\ (\operatorname{mod}\ 3k+1)$ corresponding 
    to $s=k-1$ and $s=0$. Assume
	without loss of generality that $j$ precedes $\ell$ in clockwise order around
	$C_{3k+1}$. (Thus, for the subsets $\{0,2\}$ and $\{11,0\}$ of $C_{13}$, for
	example, $0$ precedes $2$, while $11$ precedes $0$.) With arithmetic performed
	modulo $3k+1$, the neighbours of $\{j,\ell\}$ are $\{j,\ell+3\}$ and
	$\{j-3,\ell\}$, and $\{j,\ell\}$ has degree $2$ in $\mathfrak{B}_{k}$. (Thus,
	$\{0,2\}$ is adjacent to $\{0,5\}$ and $\{10,2\}$ in $\mathfrak{B}_{4}$, while
	$\{11,0\}$ is adjacent to $\{11,3\}$ and $\{8,0\}$.)
	
	\item Suppose $j-\ell\not \equiv 2$ or $-2\ (\operatorname{mod}\ 3k+1)$. Then
	the neighbours of $\{j,\ell\}$ in $\mathfrak{B}_{k}$ are $\{j-3,\ell
	\},\{j+3,\ell\},\{j,\ell-3\}$ and $\{j,\ell+3\}$ (arithmetic modulo $3k+1$),
	and $\{j,\ell\}$ has degree $4$ in $\mathfrak{B}_{k}$. (For example, the
	neighbours, in $\mathfrak{B}_{4}$, of $\{0,5\}$ are $\{10,5\},\{3,5\},\{0,2\}$
	and $\{0,8\}$).
\end{enumerate}

We emphasize that the sets $\bk{j,\ell}$ are unordered.  For example, the neighbours of $\bk{2,7}$ in $\bcl_3$ are $\bk{2,0} = \bk{0,2}, \bk{2,4},\bk{9,7} = \bk{7,9}$ and $\bk{5,7}$.  Examples of these graphs are given in Figures \ref{fig:i:iC7},	\ref{fig:i:iC10},  \ref{fig:i:iC13},  and \ref{fig:i:iC16} below, but where the set braces are removed to reduce visual clutter.  			

The two-number identifiers on the vertices of a bracelet graph and the $i$-sets of a cycle $C_{3k+1}$ are no coincidence; in the following series of lemmas and observations, we show that they are one and the same.  That is, we show that the vertex $\{i,j\}$ in the bracelet graph $\bcl_k$ corresponds to the $i$-set of $C_{3k+1}$ represented by $\cid{i,j}$, and hence that $\ig{C_{3k+1}} \cong \bcl_k$.


\begin{figure}[H] \centering	 	 
	\begin{tikzpicture}	[scale=.8]		
		
		\foreach \c / \i / \j in {
			0/	0	/	2	,
			1/	3	/	5	,
			2/	6	/	1	,
			3/	2	/	4		
		}
		{
			\node[std,label={0: \scriptsize ${\i,\j}$}] (v\c) at (2*\c, 2) {}; 					
		}				
		
		\foreach \c / \i / \j in {
			0/	0	/	5	,
			1/	3	/	1	,
			2/	6	/	4								
		}
		{
			\node[std,label={0: \scriptsize ${\i,\j}$}] (z\c) at (1+2*\c, 0) {}; 					
		}			
		
		\foreach \i / \j in {0/1,1/2,2/3}
		{
			\draw (z\i)--(v\j);
			\draw (z\i)--(v\i);
		}
		
		\draw (v0) to [out=270,in=160]  (1,-1) to [out=0, in=180]  (5,-1)  to [out=10,in=270]  (v3); 
		
	\end{tikzpicture} 
	\caption{$\bcl_2 \cong \ig{C_{7}}$.}
	\label{fig:i:iC7}		 
\end{figure}

\begin{figure}[H] \centering	 	 
	\begin{tikzpicture}	[scale=.8]		\hspace{3mm}
		
		\foreach \c / \i / \j in {
			0/	7	/	9	,
			1/	0	/	2	,
			2/	3	/	5	,
			3/	6	/	8	,
			4/	1	/	9	}
		{
			\node[std,label={0: \scriptsize ${\i,\j}$}] (v\c) at (2*\c, 2) {}; 					
		}				
		
		\foreach \c / \i / \j in {
			0/	2	/	7	,
			1/	0	/	5	,
			2/	3	/	8	,
			3/	1	/	6	,
			4/	4	/	9									
		}
		{
			\node[std,label={0: \scriptsize ${\i,\j}$}] (w\c) at (1+2*\c, 1) {}; 					
		}	
		
		\foreach \c / \i / \j in {
			0/	2	/	4	,
			1/	5	/	7	,
			2/	0	/	8	,
			3/	1	/	3	,
			4/	4	/	6							
		}
		{
			\node[std,label={0: \scriptsize ${\i,\j}$}] (z\c) at (2*\c, 0) {}; 					
		}			
		
		\draw	(v0)--(w0)--(z0);
		\draw	(v1)--(w1)--(z1);
		\draw	(v2)--(w2)--(z2);
		\draw	(v3)--(w3)--(z3);
		\draw	(v4)--(w4)--(z4);
		
		\draw	(v1)--(w0)--(z1);
		\draw	(v2)--(w1)--(z2);
		\draw	(v3)--(w2)--(z3);
		\draw	(v4)--(w3)--(z4);
		
		\draw (v0) to [out=55,in=45]  (w4); 
		\draw (z0) to [out=-55,in=-45]  (w4); 
		
	\end{tikzpicture} 
	\caption{$\bcl_3 \cong \ig{C_{10}}$.}
	\label{fig:i:iC10}		 
\end{figure}

\begin{figure}[H]  \begin{minipage}{0.1\textwidth} 
	\begin{tikzpicture}	[scale=.8]	
		
		\foreach \c / \i / \j in {
			0/	0	/	2	,
			1/	3	/	5	,
			2/	6	/	8	,
			3/	9	/	11	,
			4/	1	/	12	,
			5/	2	/	4	
		}
		{
			\node[std,label={0: \scriptsize ${\i,\j}$}] (t\c) at (1+2*\c, 3) {}; 					
		}				
		
		\foreach \c / \i / \j in {
			0/	2	/	10	,
			1/	0	/	5	,
			2/	3	/	8	,
			3/	6	/	11	,
			4/	1	/	9	,
			5/	4	/	12	,
			6/	2	/	7										
		}
		{
			\node[std,label={0: \scriptsize ${\i,\j}$}] (v\c) at (2*\c, 2) {}; 					
		}	
		
		\foreach \c / \i / \j in {
			0/	5	/	10	,
			1/	0	/	8	,
			2/	3	/	11	,
			3/	1	/	6	,
			4/	4	/	9	,
			5/	7	/	12									
		}
		{
			\node[std,label={0: \scriptsize ${\i,\j}$}] (w\c) at (1+2*\c, 1) {}; 					
		}					
		
		\foreach \c / \i / \j in {
			0/	5	/	7	,
			1/	8	/	10	,
			2/	0	/	11	,
			3/	1	/	3	,
			4/	4	/	6	,
			5/	7	/	9	,
			6/	10	/	12										
		}
		{
			\node[std,label={0: \scriptsize ${\i,\j}$}] (z\c) at (2*\c, 0) {}; 					
		}			
		
		
        \draw (v0)--(w0)--(v1)--(w1)--(v2)--(w2)--(v3)--(w3);
        \draw (w3)--(v4)--(w4)--(v5)--(w5)--(v6);
		
		\draw (v0) to [out=135,in=45]  (v6); 
		\draw[very thick,red,dashed] (v0) to [out=-135,in=160]  (0.2,-1) to [out=-20,in=210] (z6); 
		\draw[very thick,red,dashed] (v6) to [out=-25,in=30]  (13,-0.5) to [out=210,in=-30] (z0);

		\foreach \i / \j in {0/1,1/2,2/3,3/4,4/5,5/6}
		{
			\draw[very thick,red,dashed] (t\i)--(v\j);
			\draw[very thick,red,dashed] (t\i)--(v\i);
		}		
		
		\foreach \i / \j in {0/1,1/2,2/3,3/4,4/5,5/6}
		{
			\draw[very thick,red,dashed] (w\i)--(z\j);
			\draw[very thick,red,dashed] (z\i)--(w\i);
		}

	\end{tikzpicture} 
\end{minipage}
	\caption{$\bcl_4 \cong \ig{C_{13}}$ with Hamiltonian cycle in red.}
	\label{fig:i:iC13}
\end{figure}

\begin{figure}[H] 	 \centering
	\begin{tikzpicture}	[scale=0.8]	 
		
		\foreach \c / \i / \j in {
			0/	0	/	2	,
			1/	3	/	5	,
			2/	6	/	8	,
			3/	9	/	11	,
			4/	12	/	14	,
			5/	1	/	15	,
			6/	2	/	4	,
			7/	5	/	7			
		}
		{
			\node[std,label={0: \scriptsize ${\i,\j}$}] (r\c) at (2*\c, 4) {}; 					
		}			
		
		\foreach \c / \i / \j in {
			0/	0	/	5	,
			1/	3	/	8	,
			2/	6	/	11	,
			3/	9	/	14	,
			4/	1	/	12	,
			5/	4	/	15	,
			6/	2	/	7	,
			7/	5	/	10				
		}
		{
			\node[std,label={0: \scriptsize ${\i,\j}$}] (t\c) at (1+2*\c, 3) {}; 					
		}		
		
		\foreach \c / \i / \j in {
			0/	5	/	13	,
			1/	0	/	8	,
			2/	3	/	11	,
			3/	6	/	14	,
			4/	1	/	9	,
			5/	4	/	12	,
			6/	7	/	15	,
			7/	2	/	10										
		}
		{
			\node[std,label={0: \scriptsize ${\i,\j}$}] (v\c) at (2*\c, 2) {}; 					
		}	
		
		\foreach \c / \i / \j in {
			0/	8	/	13	,
			1/	0	/	11	,
			2/	3	/	14	,
			3/	1	/	6	,
			4/	4	/	9	,
			5/	7	/	12	,
			6/	10	/	15	,
			7/	2	/	13												
		}
		{
			\node[std,label={0: \scriptsize ${\i,\j}$}] (w\c) at (1+2*\c, 1) {}; 					
		}					
		
		\foreach \c / \i / \j in {
			0/	8	/	10	,
			1/	11	/	13	,
			2/	0	/	14	,
			3/	1	/	3	,
			4/	4	/	6	,
			5/	7	/	9	,
			6/	10	/	12	,
			7/	13	/	15										
		}
		{
			\node[std,label={0: \scriptsize ${\i,\j}$}] (z\c) at (2*\c, 0) {}; 					
		}			
		
		\draw	(r0)--(t0)--(v0)--(w0)--(z0);
		\draw	(r1)--(t1)--(v1)--(w1)--(z1);
		\draw	(r2)--(t2)--(v2)--(w2)--(z2);
		\draw	(r3)--(t3)--(v3)--(w3)--(z3);
		\draw	(r4)--(t4)--(v4)--(w4)--(z4);
		\draw	(r5)--(t5)--(v5)--(w5)--(z5);
		\draw	(r6)--(t6)--(v6)--(w6)--(z6);
		\draw	(r7)--(t7)--(v7)--(w7)--(z7);
		
		\draw	(r1)--(t0)--(v1)--(w0)--(z1);
		\draw	(r2)--(t1)--(v2)--(w1)--(z2);
		\draw	(r3)--(t2)--(v3)--(w2)--(z3);
		\draw	(r4)--(t3)--(v4)--(w3)--(z4);
		\draw	(r5)--(t4)--(v5)--(w4)--(z5);
		\draw	(r6)--(t5)--(v6)--(w5)--(z6);
		\draw	(r7)--(t6)--(v7)--(w6)--(z7);

		\draw (v0) to [out=155,in=190]  (0.8,5) to [out=10, in=180]  (8,5.5)  to [out=0,in=35]  (t7);    
		\draw[red] (v0) to [out=-155,in=-190]  (0.8,-1) to [out=-10, in=180]  (8,-1.5)  to [out=0,in=-35]  (w7);

		\draw[color=blue] (r0) to [out=70,in=190]  (3,6) to [out=10, in=160]  (15,5.5)  to [out=-20,in=25]  (w7);    
		\draw[color=ForestGreen] (z0) to [out=-70,in=170]  (3,-2) to [out=-10, in=200]  (15,-1.5)  to [out=20,in=-25]  (t7);    		
		
	\end{tikzpicture} 
	\caption{$\bcl_5 \cong \ig{C_{16}}$.}
	\label{fig:i:iC16}	
\end{figure}


\begin{lemma} \label{lem:i:bcl}
	For $k \geq 1$,		
		$	\ig{C_{3k+1}} \cong 	\bcl_{k}$.
\end{lemma}

To prove Lemma \ref{lem:i:bcl}, we first establish several lemmas, beginning with a formulation of the observation from page \pageref{sec:i:cycles}, for the case $n=3k+1$.  Going forward, we assume that all arithmetic in the $i$-set notation is performed modulo $3k+1$.

\begin{obs} \label{obs:i:2ver}
	Each $i$-set $X$ of $C_{3k+1}$ contains exactly two pairs of vertices, say $v_{j-1},v_{j+1}$ and $v_{\ell-1},v_{\ell+1}$,
	such that each pair has a common neighbour in $C_{3k+1}-X$; all other pairs of consecutive
	vertices of $X$ on $C_{3k+1}$ are separated by two vertices of $C_{3k+1}-X$. 
\end{obs}

\begin{lemma} \label{lem:i:Fdegree}
	Each vertex of $\ig{C_{3k+1}}$ has degree $2$ or $4$.  
\end{lemma}

\begin{proof}
	For each $i$-set $X$ of $C_{3k+1}$, we  deduce from Observation \ref{obs:i:2ver} that $X$ is one of two types (with notation as in Observation \ref{obs:i:2ver}):
	
	\begin{enumerate}[label=]
		\item \textbf{Type 1:}\ $\{v_{j-1},v_{j+1}\}\cap\{v_{\ell-1},v_{\ell+1}\}\neq\varnothing$; in this case
		$|\{v_{j-1},v_{j+1}\}\cap\{v_{\ell-1},v_{\ell+1}\}|=1$. Say ${j+1}={\ell-1}$. In our wide-angled notation, Type 1
		$i$-sets are of the form $X=\left\langle j,j+2 \right\rangle$. 
		Then the subpath $(v_{j-1},v_j,v_{j+1},v_{j+2},v_{j+3})$  of $C_{3k+1}$ has tokens on $v_{j-1}, v_{j+1},$ and $v_{j+3}$. A
		token on $v_{j-1}$ can only slide counterclockwise to $v_{j-2}$, a token on $v_{j+3}$ can only slide clockwise to $v_{j+4}$, and a token on $v_{j+1}$ is frozen. Hence $X$ has degree $2$ in $\ig{C_{3k+1}}$ (all other tokens are frozen).
		
		\item \textbf{Type 2:}\ $\{v_{j-1},v_{j+1}\}\cap\{v_{\ell-1},v_{\ell+1}\}=\varnothing$. When $X$ is a Type 2
		$i$-set, and $v_{j-1},v_{j+1},v_{\ell-1},$ and $v_{\ell+1}$ occur in this order in a clockwise direction on
		$C_{3k+1}$, then each of $v_{j-1}$ and $v_{\ell-1}$ is immediately preceded
		(counterclockwise) by two vertices of $C_{3k+1}-X$, and each of $v_{j+1}$ and $v_{\ell+1}$ is
		immediately followed (clockwise) by two vertices of $C_{3k+1}-X$. Hence tokens
		on $v_{j-1}$ and $v_{\ell-1}$ can slide counterclockwise to $v_{j-2}$
		and $v_{\ell-1}$, respectively, while tokens on $v_{j+1}$ and
		$v_{\ell+1}$ can slide clockwise to $v_{j+2}$ and
		$v_{\ell+2}$, respectively. Hence $X$ has degree $4$ in
		$\ig{C_{3k+1}}$.
	\end{enumerate}
\end{proof}


The following lemma is straightforward from the above proof; its proof can be found in \cite[Lemma 4.11]{LauraD}.

\begin{lemma} \label{lem:i:bkNeigh}   
	
	Let $X$ be an $i$-set of $C_{3k+1}$.	
	\begin{enumerate}[label=(\roman*)]
		\item  \label{lem:i:bkNeigh:i}    If $X=\left\langle j,j+2\right\rangle $ for some $j\in
		\{0,...,3k\}$, then the neighbours of $X$ in $\ig{C_{3k+1}}$ are
		$\left\langle j,j+5\right\rangle $ and $\left\langle j-3,j+2\right\rangle $. 
		
		\item \label{lem:i:bkNeigh:ii}  If $X=\left\langle j,\ell\right\rangle $ for some $j\in
		\{0,...,3k\}$ and some $\ell$, where $j-\ell\not \equiv 2$ or
		$-2\ (\operatorname{mod}\ 3k+1)$, then $\ell=j+2+3s$, where $s\in
		\{1,...,k-2\}$. The neighbours of $X$ in $\ig{C_{3k+1}}$ are
		$\left\langle j,\ell+3\right\rangle ,\left\langle j,\ell-3\right\rangle
		,\left\langle j+3,\ell\right\rangle ,$ and $\left\langle j-3,\ell\right\rangle
		$.
	\end{enumerate}
\end{lemma}


This completes the proof of Lemma \ref{lem:i:bcl}.  Finally, combining Lemma \ref{lem:i:2Cn} and Lemma \ref{lem:i:bcl} reveals the full result for $i$-graphs of cycles.

\begin{theorem} \label{thm:i:Cn}
	For $n \geq 3$, $k \geq 0$,		
	\begin{align*}
		\ig{C_n} = 	
		\begin{cases}
			K_3 & \text{ if } n=3 \\
			3K_1 & \text{ if } n=3k \geq 6 \\
			\bcl_{k} & \text{ if } n = 3k+1\\
			C_n & \text{ if } n \equiv 2 \Mod{3}.
		\end{cases}
	\end{align*}
\end{theorem}

\subsection{Hamiltonicity of $\ig{C_n}$}	\label{sec:i:hamCycle}

In some of the figures presented in Section \ref{sec:i:cycles}, a Hamiltonian cycle or path is easily found; in others, it is not.  This leads to the problem of determining the values of $n$ for which $\ig{C_n}$ is Hamiltonian or Hamilton traceable (i.e. has a Hamiltonian path).  In most cases, this is not too difficult to determine, as we show next.  

\begin{theorem} \label{thm:i:Ham}  For $n\geq 3$,
	\begin{enumerate}[label=$\bullet$]
		\item If $n \equiv 0 \Mod{3}$ and $n \neq 3$, then $\ig{C_n}$ is disconnected.
		\item If $n \equiv 2 \Mod{3}$ or $n=3$, then $\ig{C_n}$ is trivially Hamiltonian.
		\item If $n  \equiv 4 \Mod{6}$, then $\ig{C_n}$ is neither Hamiltonian nor Hamilton traceable.
	\end{enumerate}
\end{theorem}

\begin{proof}	
	The first two cases are trivial, and so assume that $n \equiv 4 \Mod{6}$. Since $\ig{C_4} = \overline{K_2}$, it is non-Hamiltonian;
	hence, assume $n >4$, say $n = 6k+4$, $k \geq 1$.   With notation as above, we first count the
	number of vertices $\cid{j, \ell}$ of $\ig{C_n}$ such that $j \equiv \ell \Mod{2}$.  For $j=0$, the set of these
	vertices is
	\begin{align*}
		\mathscr{X}_0 = \{\cid{0,2}, \cid{0,8}, \dots, \cid{0,6k+2}\}
	\end{align*}	
	and $| \mathscr{X}_0| = k+1$.   Similarly, 	
	\begin{align*}
		\mathscr{X}_j = \{\cid{j,j+2}, \cid{j,j+8}, \dots, \cid{j,j+6k+2}\}
	\end{align*}	
	and
	\begin{align*}
		\sum^{6k+3}_{j=0} |\mathscr{X}_j| = (k+1)(6k+4).
	\end{align*}	
	But each vertex $\cid{j,\ell}, j \equiv \ell \Mod{2}$  occurs in exactly two sets, namely $\mathscr{X}_j$ and $\mathscr{X_{\ell}}$.	
	Hence,
	\begin{align*}
		\left| \bigcup_{j=0}^{6k+3} \mathscr{X}_j \right| = (k+1)(3k+2).
	\end{align*}
	
	\noindent Since 
	\begin{align*}
		|V(\ig{C_n})| = \frac{n(n-1)}{6} = (3k+2)(2k+1),
	\end{align*}
	$\ig{C_n}$ has $k(3k+2)$ vertices $\cid{j,\ell}$ such that $j \not\equiv \ell \Mod{2}$.  

	We show next that each vertex $\cid{j,\ell}$ such that $j \equiv \ell \Mod{2}$ is adjacent only to vertices $\cid{j',\ell'}$ such that $j' \not\equiv \ell' \Mod{2}$, and vice versa.

	Let $\cid{j,\ell}$ be any vertex of $\ig{C_n}$ such that $j \equiv \ell \Mod{2}$.  Then, with arithmetic in
	the subscripts performed modulo $n$, 
	\begin{align*}
		N(\cid{j,\ell}) = \begin{cases}
			\{  \cid{j, \ell+3},  \cid{j, \ell-3},  \cid{j+3, \ell},  \cid{j-3, \ell}  \} & \text{ if } \ell -j \not\equiv \pm 2 \Mod{n} \\
			\{  \cid{j, \ell+3},   \cid{j-3, \ell}  \} & \text{ if } \ell -j \equiv  2 \Mod{n} \\
			\{    \cid{j, \ell-3},  \cid{j+3, \ell} \} & \text{ if } \ell -j \equiv -2 \Mod{n}, \\			
		\end{cases}
	\end{align*}

	\noindent	where, since $n$  is even, $j \not\equiv \ell \pm 3  \Mod{2}$ and  $\ell \not\equiv j \pm 3 \Mod{2}$. Hence each vertex $\cid{j,\ell}$ of
	$\ig{C_n}$ such that $j \equiv \ell \Mod{2}$  is adjacent only to vertices  $\cid{j',\ell'}$ such that $j' \not\equiv \ell' \Mod{2}$.

	Therefore, $\ig{C_n}$ is bipartite with $(k+1)(3k+2)$ vertices in one partite set, and $k(3k+2)$ in the other.   Since the cardinalities of the partite sets differ by more than one, the result follows.  
\end{proof}

\medskip

The case for $n \equiv 1 \Mod{6}$ is more complicated.   From Figures \ref{fig:i:iC7} and \ref{fig:i:iC13} given above, $\ig{C_7}$ and $\ig{C_{13}}$ are Hamiltonian:  $\ig{C_7}$ trivially so,  and for  $\ig{C_{13}}$, illustrated in red (dashed).    For $n \equiv 1 \Mod{6}$ with $n \geq 19$, we claim that $\ig{C_n}$ is not Hamiltonian.  Consider $\ig{C_{19}}$ given in Figure \ref{fig:i:iC19} below.   Any Hamiltonian cycle on $\ig{C_{19}}$ would include all of the vertices of degree 2 and their  degree 4 neighbours, as highlighted in red (dashed) in Figure \ref{fig:i:iC19}.  However, this (proper) subset of vertices induces a cycle in $\ig{C_{19}}$, and so $\ig{C_{19}}$ is not Hamiltonian.  A similar argument follows for larger $n$  with $n \equiv 1 \Mod{6}$, as we show next.


\begin{figure}[H] \centering	 	 
	\begin{tikzpicture}[scale=0.8]			
		
		\foreach \c / \i / \j in {
			0/	0	/	2	,
			1/	3	/	5	,
			2/	6	/	8	,
			3/	9	/	11	,
			4/	12	/	14	,
			5/	15	/	17	,
			6/	18	/	1	,
			7/	2	/	4	,
			8/	5	/	7	,
			9/	8	/	10	,
			10/	11	/	13	,
			11/	14	/	16	
		}
		{
			\node[std,label={0: \scriptsize ${\i,\j}$}] (v\c) at (1.5*\c, 5) {}; 					
		}				
		
		\foreach \c / \i / \j in {
			0/	0	/	5	,
			1/	3	/	8	,
			2/	6	/	11	,
			3/	9	/	14	,
			4/	12	/	17	,
			5/	15	/	1	,
			6/	18	/	4	,
			7/	2	/	7	,
			8/	5	/	10	,
			9/	8	/	13	,
			10/	11	/	16										
		}
		{
			\node[std,label={0: \scriptsize ${\i,\j}$}] (w\c) at (0.75+1.5*\c, 4) {}; 					
		}			
		
		\foreach \c / \i / \j in {
			0/	0	/	8	,
			1/	3	/	11	,
			2/	6	/	14	,
			3/	9	/	17	,
			4/	12	/	1	,
			5/	15	/	4	,
			6/	18	/	7	,
			7/	2	/	10	,
			8/	5	/	13	,
			9/	8	/	16							
		}
		{
			\node[std,label={0: \scriptsize ${\i,\j}$}] (t\c) at (1.5+1.5*\c, 3) {}; 					
		}

		\foreach \c / \i / \j in {
			0/	0	/	11	,
			1/	3	/	14	,
			2/	6	/	17	,
			3/	9	/	1	,
			4/	12	/	4	,
			5/	15	/	7	,
			6/	18	/	10	,
			7/	2	/	13	,
			8/	5	/	16			
		}
		{
			\node[std,label={0: \scriptsize ${\i,\j}$}] (y\c) at (2.25+1.5*\c, 2) {}; 					
		}

		\foreach \c / \i / \j in {
			0/	0	/	14	,
			1/	3	/	17	,
			2/	6	/	1	,
			3/	9	/	4	,
			4/	12	/	7	,
			5/	15	/	10	,
			6/	18	/	13	,
			7/	2	/	16			
		}
		{
			\node[std,label={0: \scriptsize ${\i,\j}$}] (x\c) at (3+1.5*\c, 1) {}; 					
		}

		\foreach \c / \i / \j in {
			0/	0	/	17	,
			1/	3	/	1	,
			2/	6	/	4	,
			3/	9	/	7	,
			4/	12	/	10	,
			5/	15	/	13	,
			6/	18	/	16							
		}
		{
			\node[std,label={0: \scriptsize ${\i,\j}$}] (z\c) at (3.75+1.5*\c, 0) {}; 		
		}			
		
		
		\foreach \c in {0,1,2,3,4,5,6}
		{\draw	(w\c)--(t\c) (y\c)--(x\c);
		}
		
		\draw (w7)--(t7) (y7)--(x7);
		\draw (w8)--(t8) (y8);
		\draw (w9)--(t9);
		
		\draw (w10)--(t9) (y8)--(x7);
		\draw (w9)--(t8) (y7)--(x6);
		\draw (w8)--(t7) (y6)--(x5);
		\draw (w7)--(t6) (y5)--(x4);
		\draw (w6)--(t5) (y4)--(x3);	
		\draw (w5)--(t4) (y3)--(x2);
		\draw (w4)--(t3) (y2)--(x1);
		\draw (w3)--(t2) (y1)--(x0);
		\draw (w2)--(t1) (y0);
		\draw (w1)--(t0);
		
		

		\draw[very thick,red,dashed] (v0) to [out=250,in=160]  (3,-1) to [out=0, in=180]  (13.5,-1)  to [out=10,in=-35]  (x7);    		
		
		\draw[very thick,red,dashed]  (v11) to [out=300,in=20]  (13.5,-2.5) to [out=180, in=0]  (4.5,-2.5)  to [out=170,in=215]  (x0);    
		
		\foreach \i / \j in {0/1,1/2,2/3,3/4,4/5,5/6,6/7,7/8,8/9,9/10,10/11}
		{
			\draw[very thick,red,dashed] (w\i)--(v\j);
			\draw[very thick,red,dashed] (v\i)--(w\i);
		}

		\foreach \i / \j in {0/1,1/2,2/3,3/4,4/5,5/6,6/7}
		{
			\draw[very thick,red,dashed] (z\i)--(x\j);
			\draw[very thick,red,dashed] (z\i)--(x\i);
		}

		\draw[thick,blue,densely dotted]  (t0) to [out=250,in=160]  (3,-1.5) to [out=0, in=180]  (13.5,-1.5)  to [out=10,in=-35]  (t9); 
		\draw (w0) to [out=250,in=160]  (3,-1.7) to [out=0, in=180]  (13.5,-1.7)  to [out=10,in=-35]  (y8); 
		\draw (y0) to [out=250,in=160]  (3,-1.9) to [out=0, in=180]  (13.5,-1.9)  to [out=10,in=-35]  (w10);   
		
		\foreach \i / \j in {0/1,1/2,2/3,3/4,4/5,5/6,6/7,7/8,8/9}
		{
			\draw[thick,blue,densely dotted] (y\i)--(t\j);
			\draw[thick,blue,densely dotted] (y\i)--(t\i);
		}			
		
		
	\end{tikzpicture} 
	\caption{$\ig{C_{19}}$.}
	\label{fig:i:iC19}		 
\end{figure}


\begin{theorem} \label{thm:i:1m6nonHam}
	If $n\geq 19$ and $n \equiv 1 \Mod{6}$, then $\ig{C_n}$ is not Hamiltonian.
\end{theorem}


\begin{proof}
	As shown in Lemma \ref{lem:i:bkNeigh} \ref{lem:i:bkNeigh:i}, the neighbours of the $i$-set
	$X=\cid{0,2}$ in $\ig{C_{6k+1}}$ are
	$\cid{-3,2}$ and $\cid{ 0,5}$. On the
	other hand, the neighbours of $\cid{0,5}$ are
	$\cid{0,2}$, $\cid{0,8}$,
	$\cid{ 3,5}$ and $\cid{-3,5}$. Note
	that $\cid{0,2}$ and $\cid{3,5}$ are
	Type 1 vertices while $\cid{0,8}$ and $\cid{
		-3,5}$ are Type 2. Similarly, $\cid{-3,2}$
	has Type 1 neighbours $\cid{0,2}$ and $\cid{
		-3,-1}$, and Type 2 neighbours $\cid{ -3,5}$
	and $\cid{-6,2}$. A similar remark holds for any $i$-set
	$\cid{j,j+5},\ j\in\{0,...,3k\}$. We call these $i$-sets
	\emph{Type 2a} $i$-sets.
	
	Consider the $i$-set $X=\left\langle 0,2+3s\right\rangle $, where
	$s\in\{2,...,2k-3\}$ (hence $k\geq 3$). As proved in Lemma  \ref{lem:i:bkNeigh} \ref{lem:i:bkNeigh:ii}, the neighbours of $X$ in $\ig{C_{6k+1}}$ are the Type 2
	$i$-sets
	\[
	\cid{ 0,2+3(s+1)},\cid{ 0,2+3(s-1)}
	,\cid{ 3,2+3s},\cid{ -3,2+3s}.
	\]
	Observe that for the given range of $s$, these are all Type 2 $i$-sets. A
	similar remark holds for any $i$-set $\cid{ j,j+2+3s}$,
	where $s\in\{2,...,2k-3\}$. We call these $i$-sets \emph{Type 2b} $i$-sets.
	
	Therefore there are exactly $6k+1$ Type 1 $i$-sets of $C_{6k+1}$ and the same
	number of Type 2a $i$-sets, and each Type 1 $i$-set has exactly two neighbours
	in $\ig{C_{6k+1}}$, both of which are Type 2a $i$-sets, and,
	conversely, each Type 2a $i$-set has two Type 1 neighbours in $\ig{C_{6k+1}}$. We deduce that the subgraph of $\ig{C_{6k+1}}$ induced by
	its Type 1 vertices and their neighbours consists of only Type 1 and Type 2a
	vertices, and is $2$-regular. Hence if $\ig{C_{6k+1}}$ has Type 2b
	vertices, that is, if $k\ge3$, then $\ig{C_{6k+1}}$ is
	non-Hamiltonian. 
\end{proof}


\subsection{Traceability of $\ig{C_{6k+1}}$, where $k\geq3$.}  \label{subsec:i:C3k1Trac}

When $k\geq 3$ and $n=6k+1$, we have previously established that $\ig{C_{n}}$ has no Hamiltonian cycle.  We now instead prove that it has a Hamilton path.

\begin{theorem} \label{thm:i:C6r1}
	For $n=6k+1$, $k \geq 3$, $\ig{C_n}$ is Hamilton traceable.
\end{theorem}

\begin{proof}	
	Say $k\geq3$, and consider $C_{6k+1}$ and an $i$-set $\cid{j,\ell} $. The distance from $j$ to $\ell$ on $C_{6k+1}$ is the length of the shorter path, thus $d(j,\ell)\in
	\{2,5,8, \dots ,3k-1\}$.
	
	\begin{enumerate}[label=$\bullet$]
		\item If $k$ is even, say $k=2k^{\prime}$, then, in $C_{12k^{\prime}+1}$, we
		see that $d(j,\ell)\in\{2,5,8, \dots ,6k^{\prime}-1\}$. This set contains an equal
		number of even and odd integers.
		
		\item If $k$ is odd, say $k=2k^{\prime}+1$, then, in $C_{12k^{\prime}+7}$, we
		see that $d(j,\ell)\in\{2,5,8, \dots ,6k^{\prime}+2\}$. This set contains more even
		than odd integers.
	\end{enumerate}
	
	Consider, again, the subgraph of $\ig{C_{6k+1}}$ induced by its Type 1
	vertices and their neighbours, which is $2$-regular (as above) and has order
	$2(6k+1)$. Denote this graph by $\mathcal{H}_{2,5}$. 
	In Figure \ref{fig:i:iC19}, $\mathcal{H}_{2,5}$ is the subgraph induced by the vertices: 
	\begin{align*} \{\cid{0,2}, \cid{0,5}, \cid{3,5}, \cid{3,8}, \cid{6,8} \dots, \cid{2,16} \} 
	\end{align*}
	with the red (dashed) edges forming a cycle.
	\noindent We argue below that in general
	$\mathcal{H}_{2,5}$ is in fact connected; that is, $\mathcal{H}_{2,5}$ is a cycle.
	
	
	Note that, with arithmetic modulo $6k+1$,
	\begin{align*}
		\mathcal{W}_{2,5}  = 
		&\cid{ 0,2}, \cid{0,2+3}, \cid{3,2+3}, \cid{3,2+2\cdot3}, \cid{ 2\cdot3,2+2\cdot3},  \\
		& \dots, \cid{3x,2+3x}, \cid{3x, 2+3(x+1)} \dots
	\end{align*}	
	\noindent is a walk in $\mathcal{H}_{2,5}$. When does $\cid{ 0,2} $
	recur? There are two cases to consider, each having two subcases.
	
	
	\noindent\textbf{Case 1}:\textbf{\hspace{0.1in}}When $\cid{
		0,2} $ recurs for the first time, an even cycle is formed. Then
	$\cid{ 0,2} =\cid{ 3x,2+3x} $ for some
	integer $x$. Since $\cid{ 0,2} =\cid{
		2,0} $, there are two subcases.\smallskip
	
	\noindent\textbf{Case 1.1}:\hspace{0.1in}$3x\equiv2\ (\operatorname{mod}%
	\ 6k+1)$ and $2+3x\equiv0\ (\operatorname{mod}\ 6k+1)$, that is,
	$3x\equiv-2\ (\operatorname{mod}\ 6k+1)$. But then $2\equiv
	-2\ (\operatorname{mod}\ 6k+1)$, which is impossible because $k>0$.\smallskip
	
	\noindent\textbf{Case 1.2}:\hspace{0.1in}$3x\equiv0\ (\operatorname{mod}%
	\ 6k+1)$ and $2+3x\equiv2\ (\operatorname{mod}\ 6k+1)$. Since $\gcd
	(3,6k+1)=1$, $x\equiv0\ (\operatorname{mod}\ 6k+1)$. Then the first time
	$\cid{ 0,2} $ recurs on $\mathcal{W}_{2,5}$ is therefore
	when $x=6k+1$. It follows that $\mathcal{W}_{2,5}$ contains the cycle
	\begin{align*}
		\cid{0,2}, \cid{0,2+3}, \cid{3,2+3}, \cid{3, 2+2 \cdot 3}, \cid{2\cdot3,2+2\cdot3}, \dots, \cid{0,2}
	\end{align*}	
	of length $2(6k+1)=|V(\mathcal{H}_{2,5})|$. Therefore $\mathcal{H}_{2,5}\cong
	C_{12k+2}$.
	
	
	\noindent\textbf{Case 2}:\textbf{\hspace{0.1in}}When $\cid{
		0,2} $ recurs for the first time, an odd cycle is formed. Then
	$\cid{ 0,2} =\cid{ 3x,2+3(x+1)} $ for
	some integer $x$. Again there are two subcases.\smallskip
	
	\noindent\textbf{Case 2.1}:\textbf{\hspace{0.1in}}$3x\equiv
	2\ (\operatorname{mod}\ 6k+1)$ and $2+3(x+1)\equiv0\ (\operatorname{mod}%
	\ 6k+1)$. Then $7\equiv0\ (\operatorname{mod}\ 6k+1)$, which is impossible
	because $k>1$.\smallskip
	
	\noindent\textbf{Case 2.2}:\textbf{\hspace{0.1in}}$3x\equiv
	0\ (\operatorname{mod}\ 6k+1)$ and $2+3(x+1)\equiv2\ (\operatorname{mod}%
	\ 6k+1)$. This is likewise impossible.
	
	Therefore, we conclude in all cases that $\mathcal{H}_{2,5} \cong C_{12k+2}$.

	
	In general, for fixed $\ell\equiv2\ (\operatorname{mod}\ 6)$ and $2\leq
	\ell\leq3k-1$, denote the subgraph of $\ig{C_{6k+1}}$ induced by the
	$i$-sets of the form $\cid{ j,j+\ell} $ and $\cid{
		j,j+\ell+3} $, where $j\in\{0,...,6k\}$, by $\mathcal{H}_{\ell,\ell+3}$.
	
	Then	
	\begin{align*} 
		\mathcal{W}_{\ell, \ell+3} & = \cid{0,\ell}, \cid{0,\ell+3}, \cid{3, \ell+3}, \cid{3,\ell+2\cdot3}, \cid{2\cdot3, \ell+2\cdot3}, \\ 
		& \dots, \cid{3x, \ell+3x}, 
		\cid{3x, \ell+3(x+1)} \dots
	\end{align*}	
	is a walk in $\mathcal{H}_{\ell,\ell+3}$. When does $\cid{ 0,\ell} $
	recur? Again there are two cases to consider.
	
	\bigskip
	
	\noindent\textbf{Case 3}:\textbf{\hspace{0.1in}} When $\cid{ 0,\ell} $ recurs for the first time, an
	even cycle is formed. Then $\cid{ 0,\ell} =\cid{
		3x,\ell+3x} $ for some integer $x$. Since $\cid{
		0,\ell} =\cid{ \ell,0} $, there are two
	subcases.\smallskip
	
	\noindent\textbf{Case 3.1}:\textbf{\hspace{0.1in}}$3x\equiv
	\ell\ (\operatorname{mod}\ 6k+1)$ and $\ell+3x\equiv0\ (\operatorname{mod}\ 6k+1)$,
	that is, $3x\equiv-\ell\ (\operatorname{mod}\ 6k+1)$. Then $2\ell\equiv
	0\ (\operatorname{mod}\ 6k+1)$ and, since $\gcd(2,6k+1)=1$, $\ell\equiv
	0\ (\operatorname{mod}\ 6k+1)$. Since $2\leq \ell\leq3k-1$, this is
	impossible.\smallskip
	
	\noindent\textbf{Case 3.2}:\textbf{\hspace{0.1in}}$3x\equiv
	0\ (\operatorname{mod}\ 6k+1)$ and $\ell+3x\equiv \ell\ (\operatorname{mod}\ 6k+1)$,
	i.e., $x\equiv0\ (\operatorname{mod}\ 6k+1)$. Therefore the first time
	$\cid{ 0,\ell} $ recurs on $\mathcal{W}_{\ell,\ell+3}$ is when
	$x=6k+1$. It follows that $\mathcal{W}_{\ell,\ell+3}$ contains the cycle
	\begin{align*}
	\cid{ 0,\ell} ,\cid{ 0,\ell+3} ,\cid{
		3,\ell+3} ,...,\cid{ 3(x-1),\ell+3x} ,\cid{
		0,\ell}
	\end{align*}
	of length $2(6k+1)=|V(\mathcal{H}_{\ell,\ell+3})|$. Therefore $\mathcal{H}%
	_{\ell,\ell+3}\cong C_{12k+2}$.
	
	\bigskip
	
	\noindent\textbf{Case 4}:\textbf{\hspace{0.1in}}When $\cid{
		0,\ell} $ recurs for the first time, an odd cycle is formed. Then
	$\cid{ 0,\ell} =\cid{ 3x,\ell+3(x+1)} $ for
	some integer $x$.\smallskip
	
	\noindent\textbf{Case 4.1}:\textbf{\hspace{0.1in}}$3x\equiv
	\ell\ (\operatorname{mod}\ 6k+1)$ and $\ell+3(x+1)\equiv0\ (\operatorname{mod}%
	\ 6k+1)$. Then $2\ell+3\equiv0\ (\operatorname{mod}\ 6k+1)$, or $2\ell\equiv
	-3\equiv6k-2\ (\operatorname{mod}\ 6k+1)$. This implies that $\ell\equiv
	3k-1\ (\operatorname{mod}\ 6k+1)$ and the restrictions on $\ell$ show that
	$\ell=3k-1$. That is, there is exactly one value of $\ell$ for which these
	congruences hold. Moreover, $6x+3\equiv0\ (\operatorname{mod}\ 6k+1)$, i.e.,
	$2x\equiv-1\equiv6k\ (\operatorname{mod}\ 6k+1)$. Hence $x=3k$. 
	
	Therefore $\mathcal{W}_{\ell,\ell+3}$ contains the cycle%
	\begin{align*}
	\cid{ 0,\ell} ,\cid{ 0,\ell+3} ,\cid{
		3,\ell+3} ,...,\cid{ 3x,\ell+3x} ,\cid{
		\ell,0}
	\end{align*}
	of length $2x+1=6k+1$. Consider the distances $d(0,3k-1)$ and $d(0,3k+2)$ on
	$C_{6k+1}$. Observe that $d(0,\ell)=d(0,3k-1)=3k-1$ and
	$d(0,\ell+3)=d(0,3k+2)=6k+1-(3k+2)=3k-1$. It follows that $\mathcal{H}_{\ell,\ell+3}$
	consists of all $i$-sets $\cid{ p,q} $ such that, on
	$C_{6k+1}$, $d(p,q)=3k-1$, and there are exactly $6k+1$ such $i$-sets. Hence
	$|V(\mathcal{H}_{\ell,\ell+3})|=6k+1$, that is, $\mathcal{H}_{\ell,\ell+3}$ is exactly the
	cycle $\cid{ 0,\ell} ,\cid{ 0,\ell+3}
	,\cid{ 3,\ell+3} ,...,\cid{ 3x,\ell+3x}
	,\cid{ \ell,0} $. 
	
	Since $\ell=3k-1$ and $\ell\equiv2\ (\operatorname{mod}\ 6)$, we deduce that $k$ is
	odd; say $k=2k^{\prime}+1$, where $k^{\prime}\geq1$. Then $\ell=6k^{\prime}+2$.
	The smallest cycle $C_{6k+1}$ where $k\geq3$ for which this case occurs is
	$C_{19}$, in which case $\ell=8$ (see Figure \ref{fig:i:iC19}).
	
	\noindent\textbf{Case 4.2}:\textbf{\hspace{0.1in}}$3x\equiv
	0\ (\operatorname{mod}\ 6k+1)$ and $\ell+3(x+1)\equiv \ell\ (\operatorname{mod}%
	\ 6k+1)$. From the first congruence, $x\equiv0\ (\operatorname{mod}\ 6k+1)$,
	and so, from the second congruence, $3\equiv0\ (\operatorname{mod}\ 6k+1)$.
	This is impossible.
	
	\paragraph*{To summarize:}
	
	Fix $\ell\in\{2,5,8,...,3k-1\}$.  	
	\begin{enumerate}[label=$\bullet$]
		\item If $k=2k^{\prime}$, then by (1) and Cases 3 and 4, the subgraphs
		$\mathcal{H}_{2,5},...,\mathcal{H}_{6k^{\prime}-4,6k^{\prime}-1}$ of
		$\ig{C_{12k^{\prime}+1}}$ all have order $24k^{\prime}+2$, and
		$\mathcal{H}_{2,5}\cong\cdots\cong\mathcal{H}_{6k^{\prime}-4,6k^{\prime}%
			-1}\cong C_{24k^{\prime}+2}$.
		
		\item If $k=2k^{\prime}+1$, then by (2) and Cases 3 and 4, the subgraphs
		$\mathcal{H}_{2,5},...,\mathcal{H}_{6k^{\prime}-4,6k^{\prime}-1}$ of
		$\ig{C_{12k^{\prime}+7}}$ all have order $24k^{\prime}+14$, and
		$\mathcal{H}_{2,5}\cong\cdots\cong\mathcal{H}_{6k^{\prime}-4,6k^{\prime}%
			-1}\cong C_{24k^{\prime}+14}$. However, the subgraph $\mathcal{H}_{6k^{\prime
			}+2,6k^{\prime}+5}$ of $\ig{C_{12k^{\prime}+7}}$ has order
		$12k^{\prime}+7$ and $\mathcal{H}_{6k^{\prime
			}+2,6k^{\prime}+5}\cong C_{12k^{\prime}+7}$.  In Figure~\ref{fig:i:iC19},
        the subgraph $\mathcal{H}_{8,11}$ of $\ig{C_{19}}$ is shown with blue (dotted) edges. 
	\end{enumerate}	
	\noindent 	In either case, each vertex of $\ig{C_{6k+1}}$ belongs to $\mathcal{H}_{\ell,\ell+3}$
	for some $\ell\equiv2\ (\operatorname{mod}\ 6)$. 
	
	\subsubsection*{Connecting the Subgraphs $\mathcal{H}_{\ell,\ell+3}$ to Form a
		Hamilton Path of $\ig{C_{6k+1}}$}
	
	Denote the subgraph of $\ig{C_{6k+1}}$ that consists of the union of
	the cycles $\mathcal{H}_{\ell,\ell+3}$ by $\mathcal{H}$, and the set of edges of
	$\ig{C_{6k+1}}$ that do not belong to $\mathcal{H}$ by $\mathcal{E}$.
	Since each vertex of $\ig{C_{6k+1}}$ belongs to $H_{\ell,\ell+3}$ for some
	$\ell\equiv2\ (\operatorname{mod}\ 6)$, $\mathcal{H}$ is a spanning subgraph of
	$\ig{C_{6k+1}}$. We consider two cases, depending on whether $k$ is
	even or odd.

	\noindent\textbf{Case 1:\hspace{0.1in}}$k=2k^{\prime}$. Then
	\begin{align*}
	\mathcal{P}:\cid{ 0,2} ,\cid{ 0,5}
	,\cid{ 0,8} ,...,\cid{ 0,6k^{\prime}-4}
	,\cid{ 0,6k^{\prime}-1}
	\end{align*}
	is a path in $\ig{C_{12k^{\prime}+1}}$ whose edges belong alternately
	to $\mathcal{H}$ and to $\mathcal{E}$, beginning with the edge $(\cid{
		0,2} ,\cid{ 0,5} )$ in $\mathcal{H}_{2,5}$ and
	ending with the edge $(\cid{ 0,6k^{\prime}-4} ,\cid{
		0,6k^{\prime}-1} )$ in $\mathcal{H}_{6k^{\prime}-4,6k^{\prime}-1}%
	$. Moreover, $\mathcal{P}$ contains at least one vertex of each $\mathcal{H}%
	_{\ell,\ell+3}$. Let $\mathcal{T}$ be the subgraph of $\ig{C_{12k^{\prime
			}+1}}$ obtained by deleting all edges of $\mathcal{P}$ from $\mathcal{H}$,
	then adding the edges of $E(\mathcal{P})\cap\mathcal{E}$. Observe that
	$\mathcal{T}$ is a spanning subgraph of $\ig{C_{12k^{\prime}+1}}$.
	Since the edges of $\mathcal{P}$ were alternately deleted and added, all
	vertices of $\mathcal{T}$ have degree $2$, except for $\cid{
		0,2} $ and $\cid{ 0,6k^{\prime}-1} $, which
	have degree $1$. Also, by construction, $\mathcal{T}$ is connected. Therefore,
	$\mathcal{T}$ is a Hamiltonian path of $\ig{C_{12k^{\prime}+1}}$.

	\noindent\textbf{Case 2:\hspace{0.1in}}$k=2k^{\prime}+1$. The argument is similar.

	Therefore, in all cases $\ig{C_{6k+1}}, k\geq 3$ is Hamilton traceable.
\end{proof}

\medskip

This completes the characterization of which cycles have Hamiltonian or Hamiltonian
traceable $i$-graphs.

\section{Open Problems}
We conclude with a few open problems. Although the problems are stated here for $i$-graphs, many are relevant to other reconfiguration graphs pertaining to domination-type parameters and are also mentioned in \cite{MT18, LauraD}.



	


\begin{problem} \label{con:op:tree}
	\emph{Determine the structure of  $i$-graphs of various families of trees.  For example, consider 
 \begin{enumerate}
    \item[(a)] caterpillars in which every vertex has degree $1$ or $3$,
    \item[(b)] spiders ($K_{1,r}$ with each edge subdivided).
    \end{enumerate}}

\end{problem}

\begin{problem} \label{con:op:ham}
	\emph{Find more classes of $i$-graphs that are Hamiltonian, or Hamiltonian traceable.}
\end{problem}

\begin{problem} \label{con:op:domChain}
	\emph{Suppose $G_1, G_2,\dots$ are graphs such that $\ig{G_1} \cong G_2$, $\ig{G_2}\cong G_3$, $\ig{G_3} \cong G_4, \dots$ .  Under which conditions does there exist an integer $k$ such that $\ig{G_k} \cong G_1$?}
	
\end{problem}

As a special case of Problem \ref{con:op:domChain}, note that for any $n\geq 1$, $\ig{K_n} \cong K_n$, and that for $k \equiv 2 \Mod{3}$, $\ig{C_k} \cong C_k$.   

\begin{problem} 
	\emph{Characterize the graphs $G$ for which $\ig{G} \cong G$.}
\end{problem}

\medskip

\bigskip

\noindent\textbf{Acknowledgements\hspace{0.1in}}We acknowledge the support of
the Natural Sciences and Engineering Research Council of Canada (NSERC), RGPIN-2014-04760 and RGPIN-03930-2020.

\noindent Cette recherche a \'{e}t\'{e} financ\'{e}e par le Conseil de
recherches en sciences naturelles et en g\'{e}nie du Canada (CRSNG), RGPIN-2014-04760 and
RGPIN-03930-2020.
\begin{center}
\includegraphics[width=2.5cm]{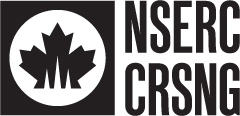}%
\end{center}
	

	\bibliographystyle{abbrv} 
	\bibliography{LTPhDBib} 
	
\end{document}